\documentclass[a4paper,11pt,english]{article}
\usepackage{amsmath}
\usepackage{amsfonts}
\usepackage{amssymb}
\usepackage{amsthm}
\usepackage[all]{xy}
\usepackage{geometry}
\usepackage{babel}

\newtheorem {thm}{Theorem}[section]
\newtheorem* {thm*}{Theorem}
\newtheorem{main}{Main Theorem}
\newtheorem {cor}[thm]{Corollary}
\newtheorem {lem}[thm]{Lemma}
\newtheorem {prop}[thm]{Proposition}

\theoremstyle{definition}
\newtheorem {rem}[thm]{Remark}
\newtheorem {defi}[thm]{Definition}

\newcommand{\id}{\mathrm{id}}

\DeclareMathOperator{\coker}{Coker}

\DeclareMathOperator{\Spec}{Spec}
\DeclareMathOperator{\End}{End}

\DeclareMathOperator{\Gal}{Gal}

\DeclareMathOperator{\Frob}{Frob}

\title{Prescribing valuations of the order of a point in the reductions of abelian varieties and tori}
\author{Antonella Perucca}
\date{}
\begin{document}
\maketitle

\begin{abstract}
Let $G$ be the product of an abelian variety and a torus defined over a number field $K$.
Let $R$ be a $K$-rational point on $G$ of infinite order. Call $n_R$ the number of connected components of the smallest algebraic $K$-subgroup of $G$ to which $R$ belongs.
We prove that $n_R$ is the greatest positive integer which divides the order of $(R \bmod\mathfrak p)$ for all but finitely many primes $\mathfrak p$ of $K$.
Furthermore, let $m>0$ be a multiple of $n_R$ and let $S$ be a finite set of rational primes.
Then there exists a positive Dirichlet density of primes $\mathfrak p$ of $K$ such that for every $\ell$ in $S$ the $\ell$-adic valuation of the order of $(R \bmod\mathfrak p)$ equals $\mathit v_\ell(m)$.
\end{abstract}

\section{Introduction}
Let $G$ be a semi-abelian variety defined over a number field $K$.
We consider reduction maps on $G$ by fixing a model for $G$ over an open subscheme of $\Spec \mathcal O$, where $\mathcal O$ is the ring of integers of $K$.

Remark that different choices of the model may affect only finitely many reductions because in fact any two models are isomorphic on a (possibly smaller) open subscheme of $\Spec \mathcal O$.


Let $R$ be a $K$-rational point on $G$. 
For all but finitely many primes $\mathfrak p$ of $K$ the reduction modulo $\mathfrak p$ is well defined on the point $R$ and the order of $(R \bmod \mathfrak{p})$ is finite. 
It is natural to ask the following question:
how does the order of $(R \bmod{\mathfrak p})$ behave if we vary $\mathfrak p$?

It is easy to see that if $R$ is non-zero then for all but finitely many primes $\mathfrak p$ of $K$ the point $(R \bmod \mathfrak{p})$ is non-zero.
A first consequence is that if $R$ is a torsion point of order $n$ then for all but finitely many primes $\mathfrak p$ of $K$ the order of $(R \bmod \mathfrak{p})$ is $n$. 
A second consequence is that if $R$ has infinite order then the order of $(R \bmod\mathfrak p)$ cannot take the same value for infinitely many primes $\mathfrak p$ of $K$.
In this paper we prove the following result:

\begin{main}\label{maintheorem}
Let $G$ be the product of an abelian variety and a torus defined over a number field $K$. Let $R$ be a $K$-rational point on $G$ of infinite order. 
Call $n_R$ the number of connected components of the smallest $K$-algebraic subgroup of $G$ containing $R$.
Then $n_R$ is the largest positive integer which divides the order of $(R \bmod\mathfrak p)$ for all but finitely many primes $\mathfrak p$ of $K$.
Furthermore, let $m>0$ be a multiple of $n_R$ and let $S$ be a finite set of rational primes. Then there exists a positive Dirichlet density of primes $\mathfrak p$ of $K$ such that for every $\ell$ in $S$ the $\ell$-adic valuation of the order of $(R \bmod\mathfrak p)$ equals $\mathit v_\ell(m)$.
\end{main}


It is interesting to see whether our result generalizes to semi-abelian varieties. In this generality we prove that for every integer $m>0$ there exists a positive Dirichlet density of primes $\mathfrak p$ of $K$ such that the order of $(R \bmod\mathfrak{p})$ is a multiple of $m$ (see Corollary~\ref{semiabelian}).
Also for all but finitely many primes $\mathfrak p$ the order of $(R \bmod\mathfrak{p})$ is a multiple of $n_R$ (see Proposition~\ref{nRdivides}).

The Main Theorem and the results in section~\ref{sectioncorollaries} (Proposition~\ref{independentcomparetorsion}, Proposition~\ref{independentexactness} and Corollary~\ref{semiabelian}) strengthen results which are in the literature:
\cite[Lemma 5]{KharePrasad}; \cite[Theorems 4.1 and 4.4]{Pink}; \cite[Theorem 3.1]{BGKdetecting} and \cite[Theorem 5.1]{Baranczuk06} in the case of abelian varieties.
Further papers concerning the order of the reductions of points are \cite{CheonHahn}, \cite{Kowalskikummer} and \cite{JonesRouse}. 

\section{Preliminaries}\label{preliminaries}
Let $G$ be a semi-abelian variety defined over a number field $K$. Let $R$ be a $K$-rational point on~$G$.
Write $G_R$ for the Zariski closure of~$\mathbb Z\cdot R$ in $G\times_K\bar{K}$ (with reduced structure).
Because $\mathbb Z\cdot R$ is dense in $G_R(\bar{K})$, it follows that $G_R$ is an algebraic subgroup of $G$ defined over $K$.
In particular for every algebraic extension $L$ of $K$ we have that $G_R$ is the smallest algebraic $L$-subgroup of $G$ such that $R$ is an $L$-rational point.
Write $G_R^0$ for the connected component of the identity of $G_R$. Then $G_R^0$ is an algebraic subgroup of $G$ defined over $K$ and $G_R^0(\bar{K})$ is divisible.
Write $n_R$ for the number of connected components of $G_R$.
The number $n_R$ does not get affected by a change of ground field: since $\mathbb Z\cdot R$ is Zariski-dense in $G_R(\bar{K})$ then every connected component of $G_R$ is a translate of $G_R^0$ by a $K$-rational point therefore it is also defined over $K$.

\begin{lem}\label{lemalgsubgroup}
Let $G$ be a semi-abelian variety defined over a number field $K$. Let $R$ be a $K$-rational point on~$G$. Then $G_{n_RR}=G_R^0$. 
Furthermore, let $H$ be a connected component of $G_R$. Then there exists a torsion point $X$ in $G_R(\bar{K})$ such that $H = X + G_R^0$.
\end{lem}
\textit{Proof.} 
Clearly $G_R^0$ contains $G_{n_RR}$. Also  $G_R^0$ is mapped to  
$G_{n_RR}$ by $[n_R]$. Because this map has finite kernel, $G_R^0$ and $G_{n_RR}$ have the same dimension. Then since $G_R^0$ is connected, we must have $G_{n_RR}=G_R^0$.

Let $P$ be any point in $H(\bar{K})$. Then $P+G_R^0=H$.
The point $n_RP$ is in $G_R^0(\bar{K})$. 
Since $G_R^0(\bar{K})$ is divisible, there exists a point $Q$ in $G_R^0(\bar{K})$ such that $n_R Q = n_R P$. Set $X = P - Q$, thus $X$ is a torsion point in $G_R(\bar{K})$. Then we have:
$$H = P + G_R^0 = P - Q + G_R^0 = X + G_R^0.$$
\hfill $\square$

\begin{prop}\label{nRdivides}
Let $G$ be a semi-abelian variety defined over a number field $K$. Let $R$ be a $K$-rational point on $G$. Then $n_R$ divides the order of $(R \bmod \mathfrak p)$ for all but finitely many primes $\mathfrak p$ of $K$.
\end{prop}
\textit{Proof.} 
Because of Lemma~\ref{lemalgsubgroup} there exist a torsion point $X$ in $G_R(\bar{K})$ and a point $P$ in $G_R^0(\bar{K})$ such that $R=P+X$.
Then clearly $n_RX$ is the least multiple of $X$ which belongs to $G_R^0(\bar{K})$.
Call $t$ the order of $X$. Let $F$ be a finite extension of $K$ where $P$ is defined and $G_R[t]$ is split. 
Fix a prime $\mathfrak p$ of $K$ and let $\mathfrak q$ be a prime of $F$ over $\mathfrak p$.
Call $m$ the order of $(R \bmod\mathfrak p)$.
Up to excluding finitely many primes $\mathfrak p$ of $K$, we may assume that the order of $(R \bmod\mathfrak q)$ is also $m$.
The equality $(mX \bmod \mathfrak q)=(-mP \bmod \mathfrak q)$
implies that $(mX \bmod \mathfrak q)$ belongs to $(G^0_R(F) \bmod \mathfrak q)$.
Then $(mX \bmod \mathfrak q)$ belongs to $(G^0_R \bmod \mathfrak q)[t]$.

Up to excluding finitely many primes $\mathfrak p$ of $K$, we may assume that the reduction modulo $\mathfrak q$ maps injectively $G_R[t]$ to $(G_R \bmod \mathfrak q)[t]$ and that it maps surjectively $G^0_R[t]$ onto $(G^0_R \bmod \mathfrak q)[t]$. See \cite[Lemma 4.4]{Kowalskikummer}.
We deduce that $mX$ belongs to $G^0_R[t]$.
Then $m$ is a multiple of $n_R$. 
This shows that for all but finitely many primes $\mathfrak p$ the order of $(R \bmod \mathfrak p)$ is a multiple of $n_R$.
\hfill $\square$

\begin{defi}\label{defindependentpoint}
Let $G$ be a semi-abelian variety defined over a number field $K$. Let $R$ be a $K$-rational point on $G$.  We say that $R$ is \textit{independent} if $R$ is non-zero and $G_R=G$.
\end{defi}
By this definition an independent point has infinite order. 
Notice that this definition does not depend on the choice of the number field $K$ such that $R$ belongs to $G(K)$.

In Remark~\ref{equivalentdefindependentpoint} we prove that if $G$ is the product of an abelian variety and a torus then $R$ is independent if and only if it is non-zero and the left $\End_K G$-module generated by $R$ is free.
Then rational points of infinite order on the multiplicative group or on a simple abelian variety are independent.

\begin{lem}\label{lemindependentpoint}
Let $G$ be a semi-abelian variety defined over a number field $K$. Let $R$ be a $K$-rational point on~$G$ of infinite order. Then the point $n_RR$ is independent in $G_R^0$.
Furthermore,  let $X$ be a torsion point in $G(K)$ and suppose that $R$ is independent. Then $R+X$ is independent.
\end{lem}
\textit{Proof.} 
By Lemma~\ref{lemalgsubgroup} we have $G_{n_RR}=G_R^0$ therefore $n_RR$ is independent in $G_R^0$.

For the second assertion, we have to prove that $G_{R+X}=G$. Call $t$ the order of $X$.
Clearly $G_{R+X}\supseteq G_{t(R+X)}=G_{tR}$. Because $G_R=G$ it suffices to show that $G_{tR}=G_R$.
Remark that $G_R$ contains $G_{tR}$ and that $G_R$ is mapped to  
$G_{tR}$ by $[t]$. Because $[t]$ has finite kernel, $G_R$ and $G_{tR}$ have the same dimension.
Because $G_R$ is connected it follows that $G_{tR}=G_R$.
\hfill $\square$

\begin{prop}\label{subgroupAxT}
Let $K$ be a number field. Let $G=A\times T$ be the product of an abelian variety and a torus defined over $K$. 
Then a connected algebraic $K$-subgroup of $G$ is the product of a $K$-abelian subvariety of $A$ and a $K$-subtorus of $T$. 
\end{prop}
\textit{Proof.} 
Let $V$ be an algebraic subgroup of $G$. Call $\pi_A$ and $\pi_T$ the projections of $V$ on $A$ and $T$ respectively.
Remark that $\pi_A(V)$ is a connected $K$-subgroup of $A$ therefore it is an abelian subvariety of $A$. Similarly $\pi_T(V)$ is a connected $K$-subgroup of $T$ therefore it is a subtorus of $T$. 
By replacing $G$ with $\pi_A(V)\times\pi_T(V)$, we may assume that $\pi_A(V)=A$ and $\pi_T(V)=T$.

Write $N_T=\pi_T(V\cap (\lbrace0\rbrace\times  T))$ and $N_A=\pi_A(V\cap (A\times\lbrace0\rbrace))$. Remark that $N_A$ and $N_T$ are $K$-algebraic subgroups of $A$ and $T$ respectively.
It suffices to show that $N_A=A$ and $N_T=T$ because in that case 
$V=A\times T$ and we are done.
To prove the assertion, we make a base change to $\bar{K}$.
Since the category of commutative algebraic $\bar{K}$-schemes is abelian (\cite[Theorem p. 315 \S 5.4 Expose $VI_A$ ]{SGA3.1}) it suffices to see that the quotients $\hat{A}=A/N_A$ and $\hat{T}=T/N_T$ are zero.
The quotient $A/N_A^0$ is an abelian variety 
 (see \cite[\S 9.5]{Polishchukbook}) and then the quotient of $A/N_A^0$ by the image of $N_A$ in $A/N_A^0$ is an abelian variety (see \cite[Theorem 4 p.72]{Mumfordbook}). Hence $\hat{A}$ is an abelian variety.
Because of \cite[Corollary \S 8.5]{Borelbook} the algebraic group $T/N_T^0$ is a torus. 
The quotient of $T/N_T^0$ by the image of $N_T$ in $T/N_T^0$ is an affine algebraic group (see \cite[Theorem \S 6.8]{Borelbook}). Hence $\hat{T}$ is an affine algebraic group.

Call $\alpha$ the composition of $\pi_A$ and the quotient map from $A$ to $\hat{A}$. Similarly call $\beta$ the composition of $\pi_T$ and the quotient map from $T$ to $\hat{T}$. The product map $\alpha\times\beta$ is a map from $V$ to $\hat{A}\times\hat{T}$.
Now we show that the projection $\pi_{\hat{A}}$ from $\alpha\times\beta(V)$ to $\hat{A}$ is an isomorphism.
Clearly $\pi_{\hat{A}}$ is an epimorphism. Since we are working in an abelian category, it suffices to show that $\pi_{\hat{A}}$ is a monomorphism. Because the map $\alpha\times\beta$ from $V$ to $\alpha\times\beta(V)$ is an epimorphism, it suffices to check that the maps $\pi_{\hat{A}}\circ (\alpha\times\beta)$ and $\alpha\times\beta$ have the same kernel.
The kernel of the first map is $V\cap (N_A\times T)$. The kernel of the second map is $V\cap (N_A\times T)\cap (A\times N_T)$. 
We show that these two group schemes are isomorphic because they have the same groups of $Z$-points for every $\bar{K}$-scheme $Z$. 
The $Z$-points of the first kernel are the pairs $(a,b)$ in $V(Z)$ such that $a$ lies in $N_A(Z)$. Since $(a,0)$ belongs to $V(Z)$ we deduce that $(0,b)$ lies in $V(Z)$ and so $b$ belongs to $N_T(Z)$. Then the two kernels have the same $Z$-points.
The proof that $\alpha\times\beta(V)$ is isomorphic to $\hat{T}$ is analogous.
We deduce that $\hat{A}$ and $\hat{T}$ are isomorphic. Since $\hat{A}$ is a complete variety while $\hat{T}$ is affine the only possible morphism from $\hat{A}$ to $\hat{T}$ is zero.
Then $\hat{A}$ and $\hat{T}$ are zero. 
\hfill $\square$\newline

For the convenience of the reader we prove the following remark. 

\begin{rem}\label{equivalentdefindependentpoint}
Let $G=A\times T$ be the product of an abelian variety and a torus defined over a number field $K$. 
Then a non-zero $K$-rational point $R$ on $G$ is independent if and only if the left $\End_K G$-module generated by $R$ is free.
\end{rem}
\textit{Proof.} 
The `only if' part is straightforward: if $\phi$ is a non-zero element of $\End_K G$ such that $\phi(R)=0$ then $\ker(\phi)$ is an algebraic subgroup of $G$ different from $G$ and containing $R$ hence containing $G_R$.
Now we prove the `if' part.
Suppose that $R$ is not independent.
Because of \cite[Proposition 1.5]{Ribet1979} the left $\End_K G$-submodule of $G(K)$ generated by $R$ is free if and only if 
the left $\End_{\bar{K}} G$-submodule of $G(\bar{K})$ generated by $R$ is free. Then to conclude we construct a non-zero element of $\End_{\bar{K}} G$ whose kernel contains the point $R$.

Clearly we may assume that $R$ has infinite order.
So $G^0_R$ is non-zero and since $R$ is not independent we have $G^0_R\neq G$.
By Proposition~\ref{subgroupAxT}, $G^0_R$ is the product of an abelian subvariety $A'$ of $A$ and a subtorus $T'$ of $T$. Then either $A'$ or $T'$ are non-zero and either $A\neq A'$ or $T\neq T'$.
If $A'$ is zero set $\phi_A=\id_A$, if $A'=A$ set $\phi_A=0$.
Otherwise by the Poincar\'e Reducibility Theorem 
there exists a non-zero abelian subvariety $B$ of $A$ such that $A'$ and $B$ have finite intersection and such that the map $$\alpha:A'\times B\rightarrow A\quad \alpha(x,y)=x+y$$ is an isogeny.
Call $d$ the degree of $\alpha$ and remark that $d$ is the order of $A'\cap B$.
Call $\hat{\alpha}$ the isogeny from $A$ to $A'\times B$ such that $\alpha\circ\hat{\alpha}=[d]$.
Call $\pi$ the projection from $A'\times B$ to $\lbrace0\rbrace\times B$.
Set $\phi_A=\alpha\circ [d]\circ\pi\circ\hat{\alpha}$. 
Remark that if $\alpha(x,y)$ is a point on $A'$ then both $x$ and $y$ are points on $A'$. Then it is immediate to see that $\phi_A$ is a non-zero element of $\End_{\bar{K}} A$ and that its kernel contains $A'$.

If $T'$ is zero set $\phi_T=\id_T$, if $T'=T$  set $\phi_T=0$. 
Otherwise, because a subtorus is a direct factor 
there exists a non-zero $\phi_T$ in $\End_{\bar{K}} T$ such that $T'$ is contained in $\ker(\phi_T)$.
Then by construction $(\phi_A\times \phi_T)\circ [n_R]$ is a non-zero element of $\End_{\bar{K}}G$ whose kernel contains $G_R$.
\hfill $\square$

\section{The method by Khare and Prasad}\label{method}

In this section we prove the following result, which will be used in section \ref{sectioncorollaries} to prove the Main Theorem. To prove this result we generalize a method by Khare and Prasad (see \cite[Lemma 5]{KharePrasad}).

\begin{thm}\label{technicalresult}
Let $G$ be the product of an abelian variety and a torus defined over a number field $K$. Let $F$ be a finite extension of $K$.
Let $R$ be an $F$-rational point on $G$ such that $G_R$ is connected.
Fix a non-zero integer $m$. 
There exists a positive Dirichlet density of primes $\mathfrak p$ of $K$ such that the following holds: there exists a prime $\mathfrak q$ of $F$ over $\mathfrak p$ such that the order of $(R \bmod{\mathfrak q})$ is coprime to $m$.
\end{thm}

Remark that if $F=K$ the theorem simply says that there exists a positive Dirichlet density of primes $\mathfrak p$ of $K$ such that the order of $(R \bmod \mathfrak p)$ is coprime to $m$.

Let $G$ be a semi-abelian variety defined over a number field $K$. For $n$ in $\mathbb N$ call $K_{\ell^n}$ the smallest extension of $K$ over which every point of $G[\ell^n]$ is defined.
Let $R$ be in $G(K)$. Then for $n$ in $\mathbb N$ call $K(\frac{1}{\ell^n}R)$ the smallest extension of $K_{\ell^n}$
over which the $\ell^n$-th roots of $R$ are defined. Clearly the extensions $K_{\ell^{n+1}}/K_{\ell^n}$ and $K(\frac{1}{\ell^n}R)/K_{\ell^n}$ are Galois.
\begin{lem}\label{degree}
Let $G$ be a semi-abelian variety defined over a number field $K$.
Let $\ell$ be a rational prime and let $n$ be a positive integer. Suppose that $G(K)$ contains $G[\ell]$. Then the degree $[K_{\ell^n}:K]$ is a power of $\ell$ and for every $R$ in $G(K)$ the degree $[K(\frac{1}{\ell^n}R):K]$ is a power of $\ell$.
\end{lem}
\textit{Proof.} 
Since the points of $G[\ell]$ are defined over $K$,
we can embed $\Gal(K_{\ell^n}/K)$ into the group of the endomorphisms of $G[\ell^n]$ fixing $G[\ell]$.
The order of this group is a power of $\ell$ since $G[\ell^n]$ is a finite abelian group whose order is a power of $\ell$.
Now we only have to prove that the degree $[K(\frac{1}{\ell^n}R):K_{\ell^n}]$ is a power of $\ell$.
We can map the Galois group of the extension $K(\frac{1}{\ell^n}R)/K_{\ell^n}$ into $G[\ell^n]$, whose order is a power of $\ell$. This is accomplished via the Kummer map
$$\phi_{n}:
\Gal(K(\frac{1}{\ell^n}R)/K_{\ell^n}) \to G[\ell^n] ;\quad
\phi_{n}(\sigma)(R)=\sigma(\frac{1}{\ell^n}R) - (\frac{1}{\ell^n}R),$$
where $\frac{1}{\ell^n}R$ is an $\ell^n$-th root of $R$.
Since two such $\ell^n$-th roots differ by a torsion point of order dividing $\ell^n$, it does not matter which root we take.
This also implies that $\phi_n$ is injective. This proves the assertion.
\hfill $\square$

\begin{lem}\label{linearlydisjoint}
Let $G$ be the product of an abelian variety and a torus defined over a number field $K$. Let $R$ be a $K$-rational point of $G$ which is independent. Then for all sufficiently large $n$ we have:
$$K(\frac{1}{\ell^n}R)\cap K_{\ell^{n+1}}=K_{\ell^n}.$$ 
\end{lem}
\textit{Proof.} 
Consider the map
$$
	\alpha_n:
	\Gal(K(\frac{1}{\ell^{n+1}}R) / K_{\ell^{n+1}})
	\to
	\Gal(K(\frac{1}{\ell^n}R) / K_{\ell^n})
$$
given by the restriction to $K(\frac{1}{\ell^n}R)$.
To prove this lemma, it suffices to show
that $\alpha_n$ is surjective for sufficiently large $n$.

It is not difficult to check that the following diagram is well defined and commutative ($\phi_n$ is the Kummer map defined in the proof of Lemma~\ref{degree}
and $\beta_n$ is induced by the diagram):
$$\xymatrix{
	0 \ar[r]
		& \Gal(K(\frac{1}{\ell^{n+1}}R)/K_{\ell^{n+1}}) \ar[r]^(.65){\phi_{n+1}} \ar[d]_{\alpha_n}
		& G[\ell^{n+1}] \ar[r] \ar[d]_{[\ell]}
		& \coker\phi_{n+1} \ar[r] \ar@{->>}[d]^{\beta_n}
		& 0
	\\
	0 \ar[r]
		& \Gal(K(\frac{1}{\ell^{n}}R)/K_{\ell^{n}}) \ar[r]^(.65){\phi_{n}}
		& G[\ell^n] \ar[r]
		& \coker\phi_{n} \ar[r]
		& 0
}$$
If $\beta_n$ is injective then $\alpha_n$ is surjective.
Since $\beta_n$ is surjective,
it suffices to prove that $\coker\phi_{n+1}$ and $\coker\phi_n$ have the same order for sufficiently large $n$. Since the order of $\coker\phi_n$ increases with $n$, it is equivalent to show that the order of $\coker\phi_n$ is bounded by a constant which does not depend on~$n$. Since we assumed that $G_R=G$, this assertion is a special case of a result by Bertrand (\cite[Theorem 1]{Bertrandgalois}). 
\hfill $\square$

\begin{lem}\label{bogomolovAxsplitT}
Let $K$ be a number field. Let $G=A\times T$ be the product of an abelian variety defined over $K$ and a torus split over $K$. Fix a rational prime $\ell$.
If $T=0$ or if $A=0$ or if $\ell$ is odd then for every sufficiently large $n>0$ there exists an element $h_\ell$ in $\Gal(\bar{K}/K)$ which acts on $G[\ell^\infty]$ via an automorphism whose set of fixed points is $G[\ell^n]$.
If $A$ and $T$ are non-zero and $\ell=2$ then for every sufficiently large $n>0$ there exists an element $h_2$ in $\Gal(\bar{K}/K)$ which acts on $G[2^\infty]$ via an automorphism whose set of fixed points is $A[2^n]\times T[2^{n+1}]$.
\end{lem}
\textit{Proof.}
If $T=0$ then the assertion is a consequence of a result by Bogomolov (\cite[Corollaire 1]{Bogomolov}).
If $A=0$, because $T$ is split over $K$ then it suffices to remark the following fact: for every sufficiently large $n>0$ the field obtained by adjoining to $K$ the $\ell^{(n+1)}$-th roots of unity is a non-trivial extension of 
the field obtained by adjoining to $K$ the $\ell^n$-th roots of unity.
Now assume that $A$ and $T$ are non-zero. 
Call $\hat{A}$ the dual abelian variety of $A$. By applying a result of Bogomolov (\cite[Corollaire 1]{Bogomolov}) to $A\times \hat{A}$ we know that if $n>0$ is sufficiently large, there exists an element $h_\ell$ in $\Gal(\bar{K}/K)$ which acts on $A\times \hat{A}[\ell^{\infty}]$ as a homothety with factor $h$ in $\mathbb Z_\ell^*$ such that  $h \equiv 1 \pmod{\ell^n}$ and $h \not\equiv 1 \pmod{\ell^{n+1}}$.
For every $n$ the Weil paring 
$$e_{\ell^n}:A[\ell^n]\times \hat{A}[\ell^n]\rightarrow \mu_{\ell^n}$$ is bilinear, non-degenerate and Galois invariant. 
Since $e_{\ell^n}$ is bilinear and non-degenerate its image contains a root of unity $\zeta$ of order $\ell^n$. Choose $X_1\in A[\ell^n]$, $X_2\in \hat{A}[\ell^n]$ such that $e_{\ell^n}(X_1,X_2)=\zeta$. By Galois invariance and bilinearity we have:
$$
\sigma(\zeta) = \sigma\bigl(e_{\ell^n}(X_1,X_2)\bigr) = e_{\ell^n}(\sigma(X_1),\sigma(X_2))
= e_{\ell^n}(h \cdot X_1, h \cdot X_2) = \zeta^{h^2}
.$$
Because $\zeta$ generates $\mu_{\ell^n}$ then $\sigma$ acts on 
$\mu_{\ell^n}$ as a homothety with factor $h^2 \pmod{\ell^{n}}$.
Clearly $h^2 \equiv 1 \pmod{\ell^{n}}$ and $h^2 \not\equiv 1 \pmod{\ell^{n+1}}$ if $\ell$ is odd.
If $\ell=2$ and $n>1$ then $h^2 \equiv 1 \pmod{2^{n+1}}$ and $h^2 \not\equiv 1 \pmod{2^{n+2}}$.
Because $T$ is split over $K$ we deduce the following: if $\ell$ is odd the set of fixed points for the automorphism of $G[\ell^\infty]$ induced by $h_\ell$ is $G[\ell^n]$; if $\ell=2$ the set of fixed points for the automorphism of $G[2^\infty]$ induced by $h_2$ is $A[2^n]\times T[2^{n+1}]$.
\hfill $\square$
\newline

\textit{Proof of Theorem \ref{technicalresult}.}
By Proposition~\ref{subgroupAxT}, $G_R$ is the product of an abelian variety $A$ and a torus $T$ defined over $F$.
Let $R'$ be a point in $G_R(\bar{F})$ such that $2R'=R$.
Since $R$ is independent in $G_R$, the point $R'$ is independent  in $G_R$.
Call $S$ the the set of the prime divisors of $m$.
Let $K'$ be a finite extension of $F$ over which $R'$ is defined, over which $T$ is split and over which $G_R[\ell]$ is split for every $\ell$ in $S$.
Apply Lemma~\ref{linearlydisjoint} to the point $R'$, the algebraic group $G_R$ and with base field $K'$.
Then for all sufficiently large $n$ and for every $\ell$ in $S$ the intersection of $K'(\frac{1}{\ell^n}R')$ and $K'_{\ell^{n+1}}$ is $K'_{\ell^n}$.
Apply Lemma~\ref{bogomolovAxsplitT} to $G_R$ with base field $K'$: we can choose $n>0$ such that the previous assertion holds and such that for every $\ell$ in $S$ there exists $h_\ell$ as in Lemma~\ref{bogomolovAxsplitT}.
Call $L$  the compositum of the fields $K'(\frac{1}{\ell^n}R')$ and the fields $K'_{\ell^{n+1}}$ where $\ell$ varies in $S$.
By Lemma~\ref{degree}, the fields $K'(\frac{1}{\ell^n}R')\cdot K'_{\ell^{n+1}}$ where $\ell$ varies in $S$ are linearly disjoint over $K'$.
Then we can construct $\sigma$ in $\Gal(L/K)$ such that for every $\ell$ in $S$ the restriction of $\sigma$ to $K'(\frac{1}{\ell^n} R')$ is the identity and such that the restriction to $K'_{\ell^{n+1}}$ of $\sigma$ and of $h_\ell$ coincide.

Let $\mathfrak p$ be a prime of $K$ which does not ramify in $L$ and such that there exists a prime $\mathfrak w$ of $L$ which is over $\mathfrak p$ and such that $\Frob_{L/K}\mathfrak w=\sigma$.
By Chebotarev's Density Theorem there exists a positive Dirichlet density of prime ideals $\mathfrak p$ of $K$ which satisfy the above conditions. Let $\mathfrak q$ be the prime of $F$ lying under $\mathfrak w$.
Fix a prime $\ell$ in $S$ and suppose that the order of $(R\bmod \mathfrak q)$ is a multiple of  $\ell$. Up to discarding finitely many primes $\mathfrak p$ the order of $(R\bmod \mathfrak w)$ is a multiple of $\ell$. 
Let $Z$ be an element of $G_R(L)$ such that $\ell^n Z=R'$. 
Then  the order of $(Z\bmod{\mathfrak w})$ is a multiple of $\ell^{n+1}$ (respectively of $\ell^{n+2}$ if $\ell=2$).
Let $a\geq1$ be such that the order of $(aZ\bmod{\mathfrak w})$ is exactly $\ell^{n+1}$ (respectively $\ell^{n+2}$ if $\ell=2$).
Up to discarding finitely many primes $\mathfrak p$ there exists a torsion point $X$ in $G_R(L)$ of order $\ell^{n+1}$ (respectively $\ell^{n+2}$ if $\ell=2$) and such that $(aZ\bmod{\mathfrak w})=(X\bmod{\mathfrak w})$. See \cite[Lemma 4.4]{Kowalskikummer}.

Up to excluding finitely many primes $\mathfrak p$, the action of the Frobenius $\Frob_{L/K}\mathfrak w$ commutes with the reduction modulo $\mathfrak w$ of $G$ hence we deduce the following: the point $(Z \bmod{\mathfrak w})$ is fixed by the Frobenius of $\mathfrak w$ while $(X \bmod{\mathfrak w})$ is not fixed. Then the point $(aZ \bmod{\mathfrak w})$ is fixed by the Frobenius of $\mathfrak w$ and we get a contradiction.
\hfill $\square$

\section{The proof of the Main Theorem and corollaries}\label{sectioncorollaries}

In this section we prove the Main Theorem and other applications of Theorem~\ref{technicalresult}.

\begin{prop}\label{independentcomparetorsion}
Let $K$ be a number field. For every $i=1,\ldots, n$ let $G_i$ be the product of an abelian variety and a torus defined over $K$ and let $R_i$ be a point in $G_i(K)$ of infinite order.
Suppose that the point $R=(R_1,\ldots,R_n)$ in $G=G_1\times\ldots\times G_n$ is such that $G_R$ is connected.
Fix a non-zero integer $m$.
For every $i=1,\ldots,n$ fix a torsion point $X_i$ in $G_i(\bar{K})$ such that the point $X=(X_1,\ldots, X_n)$ is in $G_R(\bar{K})$.
Let $F$ be a finite extension of $K$ over which $X$ is defined.
Then there exists a positive Dirichlet density of primes $\mathfrak p$
of $K$ such that the following holds: there exists a prime $\mathfrak q$ of $F$ over $\mathfrak p$ such that for every $i=1,\ldots,n$ the order of $(R_i-X_i \bmod \mathfrak q)$ is coprime to $m$.
\end{prop}
\textit{Proof.} 
By Lemma~\ref{lemindependentpoint} the point $R$ is independent in $G_R$ and the point $R'=R-X$ is independent in $G_R$. 
Since $G_{R'}=G_R$, by Proposition~\ref{subgroupAxT} the algebraic group $G_{R'}$ is the product of an abelian variety and a torus defined over $K$.
Apply Theorem~\ref{technicalresult} to $R'$ and find a positive Dirichlet density of primes $\mathfrak p$ of $K$ such that the following holds: there exists a prime $\mathfrak q$ of $F$ over $\mathfrak p$ such that the order of $(R'\bmod \mathfrak q)$ is coprime to $m$. This clearly implies the statement.
\hfill $\square$

\begin{prop}\label{independentexactness}
Let $K$ be a number field. For every $i=1,\ldots, n$ let $G_i$ be the product of an abelian variety and a torus defined over $K$ and let $R_i$ be a point in $G_i(K)$ of infinite order.
Suppose that the point $R=(R_1,\ldots,R_n)$ in $G=G_1\times\ldots\times G_n$ is independent.
Fix a finite set $S$ of rational primes. For every $i=1,\ldots,n$ fix a non-zero integer $m_i$.
Then there exists a positive Dirichlet density of primes $\mathfrak p$
of $K$ such that for every $i=1,\ldots,n$ and for every $\ell$ in $S$
the $\ell$-adic valuation of the order of $(R_i \bmod \mathfrak p)$ is $\mathit v_\ell(m_i)$.
\end{prop}
\textit{Proof.}
For every $i=1,\ldots,n$ choose a torsion point $X_i$ in $G_i(\bar{K})$ of order $m_i$ and call $X=(X_1,\ldots,X_n)$.
Let $F$ be a finite extension of $K$ over which $X$ is defined.
Call $m$ the product of the primes in $S$.
Apply Proposition~\ref{independentcomparetorsion} to $R$ and find a positive Dirichlet density of primes $\mathfrak p$ of $K$ such that the following holds: there exists a prime $\mathfrak q$ of $F$ over $\mathfrak p$ such that the order of $(R-X\bmod \mathfrak q)$ is coprime to $m$.
Fix $\mathfrak p$ as above. Up to discarding finitely many primes $\mathfrak p$, for every $i=1,\ldots,n$ the order of $(X_i\bmod \mathfrak q)$ equals $m_i$. This implies that for every $i=1,\ldots,n$ and for every $\ell$ in $S$ the $\ell$-adic valuation of the order of $(R_i\bmod \mathfrak q)$ equals $\mathit v_\ell(m_i)$.
Then up to discarding finitely many primes $\mathfrak p$, the $\ell$-adic valuation of the order of $(R_i\bmod \mathfrak p)$ equals $\mathit v_\ell(m_i)$ for every $i=1,\ldots,n$ and for every $\ell$ in $S$.
\hfill $\square$\newline

\textit{Proof of the Main Theorem.}
Call $n$ the largest positive integer which divides the order of $(R \bmod \mathfrak p)$ for all but finitely many primes $\mathfrak p$ of $K$. By Proposition~\ref{nRdivides} we know that $n_R$ divides $n$.
Now we prove that $n$ divides $n_R$.
By Lemma~\ref{lemindependentpoint}, $G_{n_RR}$ is connected hence by Proposition~\ref{subgroupAxT} it is the product of an abelian variety and a torus defined over $K$. 
Let $\ell$ be a rational prime. Apply Theorem~\ref{technicalresult} to $n_RR$ and find infinitely many primes $\mathfrak p$ of $K$ such that the $\ell$-adic valuation of the order of $(n_RR \bmod \mathfrak p)$ is $0$. Thus there exist infinitely many primes $\mathfrak p$ of $K$ such that the $\ell$-adic valuation of the order of $(R \bmod \mathfrak p)$ is less than or equal to $\mathit v_\ell(n_R)$. This shows that $n$ divides $n_R$.
Now we prove the second assertion. 

Apply Proposition~\ref{independentexactness} to $n_RR$ in $G_{n_RR}$ and find a positive density of primes $\mathfrak p$ of $K$ such that for every $\ell$ in $S$ the $\ell$-adic valuation of the order of $(n_RR \bmod \mathfrak p)$ is $\mathit v_\ell(\frac{m}{n_R})$.
Because of the first assertion, we may assume that $n_R$ divides the order of $(R \bmod \mathfrak p)$. Then for every $\ell$ in $S$ the $\ell$-adic valuation of the order of $(R \bmod \mathfrak p)$ is $\mathit v_\ell(m)$. 
\hfill $\square$

By adapting this proof straightforwardly we may remark that $n_R$ is also the largest positive integer which divides the order of $(R \bmod\mathfrak p)$ for a set of primes $\mathfrak p$ of $K$ of Dirichlet density $1$.

\begin{lem}\label{lemnew}
Let $K$ be a number field. For every $i=1,\ldots, n$ let $G_i$ be the product of an abelian variety and a torus defined over $K$. Let $H$ be an algebraic subgroup of $G_1\times\ldots\times G_n$ such that the projection $\pi_i$ from $H$ to $G_i$ is non-zero for every $i=1,\ldots,n$. Let $\ell$ be a rational prime. Then there exists $X$ in $H[\ell^\infty]$ such that $\pi_i(X)$ is non-zero for every $i=1,\ldots, n$.
\end{lem}
\textit{Proof.}
By Proposition~\ref{subgroupAxT}, up to replacing $H$ with $H^0$ we may assume that $H$ is the product of an abelian variety and a torus. For every $i=1,\ldots,n$ since the projection $\pi_i$ is non-zero, it is easy to see that there exists $Y_i$ in $H[\ell^\infty]$ such that $\pi_i(Y_i)$ is non-zero. 
The point $Y_1$ is not in the kernel of $\pi_1$. So if $n=1$ we conclude. Otherwise let $1<r\leq n$ and suppose that $\sum_{j=1}^{r-1}Y_j$ is not in the kernel of $\pi_i$ for every $i=1,\ldots,r-1$. 
Up to replacing $Y_r$ with an element in $\frac{1}{\ell^\infty}Y_r$, we may assume that for every $i=1,\ldots,r$
either $\pi_i(Y_r)$ is zero or the order of $\pi_i(Y_r)$ is greater than the order of $\pi_i(\sum_{j=1}^{r-1}Y_j)$. Then $\sum_{j=1}^{r}Y_j$ is not in the kernel of $\pi_i$ for every $i=1,\ldots,r$.
We conclude by iterating the procedure up to $r=n$.

\begin{cor}\label{semiabelian}
Let $K$ be a number field. For every $i=1,\ldots, n$ let $G_i$ be a semi-abelian variety defined over $K$ and let $R_i$ be a point on $G_i(K)$ of infinite order. Then for every integer $m>0$ there exists a positive Dirichlet density of primes $\mathfrak p$ of $K$ such that for every $i=1,\ldots,n$ the order of $(R_i \bmod \mathfrak p)$ is a multiple of $m$.
\end{cor}
\textit{Proof.}
First we prove the case where $G_i$ is the product of an abelian variety $A_i$ and a torus $T_i$ for every $i=1,\ldots,n$.
Call $S$ the set of prime divisors of $m$.
Consider the point $R=(R_1,\ldots,R_n)$ in $G=G_1\times\ldots\times G_n$. 
We may assume that $n_R=1$ by replacing $R_i$ with $n_RR_i$ and we may assume that $m$ is square-free by replacing $R_i$ with $(m/\prod_{\ell\in S} \ell) R_i$ for every $i=1,\ldots,n$.
Since $G_R$ contains $R$, the projection from $G_R$ to $G_i$ is non-zero for every $i=1,\ldots,n$ so we can apply Lemma~\ref{lemnew}. Then for every $\ell$ in $S$ there exists 
$X_\ell$ in $G_R[\ell^\infty]$ such that all the coordinates of $X_\ell$ are non-zero.
 Write $Y=\sum_{\ell\in S} X_{\ell}$. By construction $Y$ belongs to $G_R(\bar{K})_{tors}$ and for every $\ell\in S$ the order of every coordinate of $Y$ is a multiple of $\ell$.
Let $F$ be a finite extension of $K$ where $Y$ is defined.
By Proposition~\ref{independentcomparetorsion}, there exists a positive Dirichlet density of primes $\mathfrak p$ of $K$ such that the following holds: there exists a prime $\mathfrak q$ of $F$ over $\mathfrak p$ such that the order of $(R-Y \bmod \mathfrak q)$ is coprime to $m$.
Then up to discarding finitely many primes $\mathfrak p$ the order of $(R_i \bmod \mathfrak p)$ is a multiple of $\ell$ for every $\ell$ in $S$ and for every $i=1,\ldots,n$. This concludes the proof for this case.

For every $i=1,\ldots,n$ let $G_i$ be an extension of an abelian variety $A_i$ by a torus $T_i$ and call $\pi_i$ the quotient map from $G_i$ to $A_i$.
If $\pi_i(R_i)$ does not have infinite order let $R'_i$ be a non-zero multiple of $R_i$ which belongs to $T_i(K)$. If $\pi_i(R)$ has infinite order then let $R'_i=0$.
Then $(\pi_iR_i,R'_i)$ is a $K$-rational point of $A_i\times T_i$ of infinite order.
Clearly for all but finitely many primes $\mathfrak p$ of $K$ the following holds: the order of $(R_i \bmod \mathfrak p)$ is a multiple of $m$ whenever the order of $((\pi_iR_i,R'_i) \bmod \mathfrak p)$ is a multiple of $m$.
Then we reduced to the previous case.
\hfill $\square$

\section*{Acknowledgements}
Marc Hindry suggested to me that the greatest integer dividing the order of a point for almost all reductions could be related to the number of connected components of the algebraic subgroup generated by the point.
I thank Bas Edixhoven for the idea of imitating Goursat's Lemma in Proposition~\ref{subgroupAxT}.
I also thank Brian Conrad and Ren\'e Schoof for helpful discussions.


\end{document}